\pdfoutput=1
\documentclass[paper=letter, fontsize=10pt,leqno]{scrartcl}
\usepackage[body={5.7in,7.5in}]{geometry}
\usepackage[T1]{fontenc}

\usepackage[english]{babel}
\usepackage{csquotes}

\usepackage[protrusion=true,expansion=true]{microtype}				
\usepackage{amsmath}										
\usepackage{amsthm}
\usepackage{latexsym}
\usepackage{hyperref}
\usepackage{tensor} 
\usepackage[pdftex]{graphicx}
\usepackage{url}
\usepackage{amssymb}
\usepackage{bbm}
\usepackage{MnSymbol}
\usepackage{wrapfig}
 \usepackage{float}
\usepackage[pdftex]{color}
\usepackage{eucal}
\usepackage{amscd}
 \usepackage{verbatim}
 \usepackage{bigints}
 \usepackage[format=plain, font=small]{caption}
\usepackage{epstopdf}
\usepackage{comment}
\usepackage{sectsty}					  
\allsectionsfont{\centering \normalfont\scshape}	  

\usepackage{mathtools}

\usepackage{fancyhdr}
\makeatletter
\theoremstyle{definition}

\theoremstyle{plain}

\theoremstyle{plain}
\newtheorem{prop}{\protect\propositionname}
\theoremstyle{plain}

\theoremstyle{plain}
\newtheorem{thm}{\protect\theoremname}
\theoremstyle{plain}

\theoremstyle{plain}

\author{\normalfont \large 
  T. Chumley\footnote{Department of Mathematics and Statistics, Mount Holyoke College, 50 College St, South Hadley, MA 01075}, 
\ R. Feres\footnote{Department of Mathematics and Statistics, Washington University, Campus Box 1146, St. Louis, MO 63130}, 
\ L. Garcia German\footnotemark[2]
}

\makeatother

\providecommand{\assumptionname}{Assumption}
\providecommand{\corollaryname}{Corollary}
\providecommand{\definitionname}{Definition}
\providecommand{\remarkname}{Remark}
\providecommand{\lemmaname}{Lemma}
\providecommand{\propositionname}{Proposition}
\providecommand{\theoremname}{Theorem}

\begin{document}
\title{\textrm{\textmd{
Revisiting Maxwell-Smoluchowski theory:
low surface roughness in straight channels}}}
\author{{\large{}T. Chumley}\thanks{Department of Mathematics and Statistics, Mount Holyoke College, 50
College St, South Hadley, MA 01075}, {\large{}R. Feres}\thanks{Department of Mathematics, Washington University, Campus
Box 1146, St. Louis, MO 63130}, {\large{}L.~A. Garcia German}\footnotemark[2], {\large{}G. Yablonsky}\thanks{Department of Energy, Environmental  \& Chemical Engineering, McKelvey School of Engineering, Washington University, St. Louis, MO 63130.}}
\maketitle
\begin{center}
Abstract
\par\end{center}
\begin{abstract}\footnotesize
The Maxwell-Smoluchowski (MS) theory of gas diffusion is revisited here in the context of gas transport in straight channels in the Knudsen regime of large mean free path. This classical theory is based on a phenomenological model of gas-surface interaction that posits that a fraction $\vartheta$ of molecular collisions with the channel surface consists of  diffuse collisions, i.e., the direction of post-collision velocities is distributed according to the Knudsen Cosine Law, and a fraction $1-\vartheta$ undergoes specular reflection. From this assumption one obtains the value $\mathcal{D}=\frac{2-\vartheta}{\vartheta}\mathcal{D}_K$ for the self-diffusivity constant, where $\mathcal{D}_K$ is a reference value  corresponding to $\vartheta=1$. In this paper we   show that $\vartheta$ can be expressed in terms of micro- and macro-geometric parameters for a model consisting of hard spheres colliding elastically against a rigid surface with  prescribed microgeometry.

Our refinement of the MS theory is based on the observation that the classical surface scattering operator associated to the microgeometry has a canonical velocity space diffusion approximation by a generalized Legendre differential operator whose spectral theory is known explicitly. 
More specifically, starting from an explicit description of the effective channel surface microgeometry\----a concept which incorporates both the actual surface microgeometry and the molecular radius\----and using this operator approximation, we show that $\vartheta$ can be resolved into easily obtained geometric parameters, $\vartheta= \lambda h/C$, having the following interpretation: $C$ is a macroscopic parameter determined by the shape of the channel cross-section; $h$ is a  parameter that precisely captures the degree of  roughness of the effective microgeometry, and $\lambda$ is a parameter that characterizes the overall curvature  of the surface microgeometry independent of $h$. Thus $\vartheta$ is resolved as the quotient of microscopic ($\lambda h$) over  macroscopic ($C$) signature parameters of the channel geometry. The identity $\vartheta= \lambda h/C$ holds up to higher order terms in the roughness parameter $h$, so our main result better applies to well polished, or low roughness, surfaces.

 \end{abstract}

\section{Introduction}
The study of diffusion  of low-pressure gases through long channels has been a topic of scientific  interest,   for its theoretical and practical importance, since the classical work on kinetic-molecular theory by  M. Knudsen  over a century ago. Applications to porous media systems and nano scale units are at the mainstream of contemporary  engineering and technology 
in areas such as chemical, biochemical, and environmental engineering, production of batteries and semiconductors. In all such systems, a {\em channel} can be 
identified as the key unit for transport through the given media.  The transport process depends on the  structural properties of the channel surface captured by a measure of surface irregularity generally referred to as {\em surface roughness}. 

Understanding and modeling the transport process in channels taking into account surface roughness is critical for creating new technologies. Obtaining highly polished, very low roughness silicon wafers is a major concern in the development of semiconductor devices \cite{Mo2020}. 
 Regarding porous systems, the traditional experimental methods for measuring surface roughness, both through contact and non-contact methods (optical and spectroscopic methods), are not  applicable. This is why  the development of a theoretical framework for the precise characterization of  surface roughness based on exit flow measurements is a topic of current interest. This paper is a contribution towards such characterization in situations where the nano scale units of transport are straight channels whose surface roughness is relatively low.

\subsection{The classical theory}
 We recall that in the large Knudsen number limit, when the mean free path is greater than the channel diameter, gas-surface interactions predominate over collisions between gas molecules; thus the geometric  characteristics of the channel surface  become an important factor influencing the speed of transport.  This influence is often captured phenomenologically through the introduction of a {\em tangential momentum accommodation coefficient} $\vartheta$. A popular phenomenological model,  due to  Maxwell and Smoluchowski,  is based on the assumption that all collisions of gas molecules with the  channel surface  are either purely diffusive or purely specular, with a fraction $\vartheta$ of the collisions being diffuse.  We refer to $\vartheta$ in this paper as the {\em Maxwell-Smoluchowski parameter}. From this assumption one obtains the value
\begin{equation}\label{MStheta}\mathcal{D} = \frac{2-\vartheta}{\vartheta}\mathcal{D}_K \end{equation}
for the constant of self-diffusivity, in which $\mathcal{D}_K$ is the diffusivity   obtained under the assumption that, at each   collision, the
post-collision velocity is independent of the pre-collision velocity  and satisfies the so-called Knudsen cosine law distribution, given explicitly in Equation (\ref{standard}).
See \cite{ACM2003} for a detailed theoretical discussion and further elaboration on the Maxwell-Smoluchowski model  for two-dimensional channels. See also \cite{YMN2012}, for example, for a reference on how the tangential momentum accommodation coefficient is measured experimentally for various materials and gas species in long tubes.  Reference \cite{AP} provides a survey of experimental results.

\subsection{The surface scattering operator and self-diffusivity}
A more fundamental description of surface-molecule interaction can be obtained by introducing a Markov (classical scattering) operator $P$, 
which gives the probability distribution of molecular post-collision velocities conditional on the pre-collision velocity. See, for example, \cite{CF2012} for a general description of $P$\----a self-adjoint operator on an appropriate Hilbert space\----for explicit mathematical interaction models.   

All the relevant surface-molecule properties that can affect diffusivity pertain to the spectrum of  $P$.   
As shown in \cite{CFZ2016}, one has
\begin{equation}\label{DDK}{\mathcal{D}}=\left(\int_{0}^2 \frac{2-\vartheta}{\vartheta}\, d\mu(\vartheta) \right) {\mathcal{D}_K}\end{equation}
where $\mu$ is a certain measure derived from the spectral resolution of $P$.   Although conceptually useful, the practical value of this formula is not ideal since a detailed determination of $\mu$  from an explicit molecule-surface interaction model is often not easy to obtain. 

Given the usefulness of the Maxwell-Smoluchowski phenomenological  model, it is natural to ask whether it is possible to identify 
specific surface characteristics  making up the single model parameter $\vartheta$ in Equation (\ref{MStheta}).
It turns out that a useful estimation of the quantity $\eta=\mathcal{D}/\mathcal{D}_k$ (the expression in parentheses in Equation (\ref{DDK})) and further qualitative understanding of $\vartheta$ can be achieved when the surface-molecule interaction is relatively weak; that is to say, when the surface exhibits  low roughness.  Since we are here mostly concerned with geometric characteristics of the surface affecting diffusivity, we assume that collisions are elastic and no energy is exchanged between surface and molecules, in which case weak scattering should be understood in the sense that the channel surface has a relatively high degree of polish, or is fairly flat.

\subsection{Signature geometric parameters of the channel}
Such an estimation of $\eta$ requires in the first place a precise characterization of what is to be understood by {\em roughness}.
One of the main contributions of this paper is to provide a mathematical characterization of roughness  that precisely connects the details of   surface-molecule interaction  and  $\eta$ under conditions of relatively weak  interaction.   This is captured by the geometric parameter $h$ defined below. Since the term `roughness' is already widely used and our $h$ is very specific,  we will refer to it here instead as the {\em flatness parameter}. (Small values of $h$ correspond to low roughness.) The main observation of this paper is as follows. An analysis of the operator $P$ based on approximating it by a diffusion operator in velocity space, which will be detailed shortly, reveals that
a small set of parameters enters into the description of   $\vartheta$ (as defined in Equation (\ref{MStheta})). These parameters, which will be defined mathematically later in the paper, are:
\begin{itemize}
\item The {\em flatness parameter} $h$: an overall measure of how flat, or polished, the surface is. 
\item The {\em shape parameters} $\lambda_1, \lambda_2$. These are a measure of mean surface curvature not affected by 
the flatness per se, as will be better explained below. When $\lambda_1=\lambda_2$, we say that the surface microgeometry is {\em isotropic}. The assumption of isotropic microgeometry will be made throughout this paper and we write $\lambda=\lambda_1=\lambda_2$. Our main result, the determination of $\mathcal{D}$ in terms of channel geometric parameters, can still be stated in the non-isotropic case except that the self-diffusivity constant would then depend on the direction tangent to the surface in which it is measured. For simplicity of presentation we do not consider this more general case here.
\item The {\em macroscopic parameter} $C$. This parameter accounts for the shape of the channel cross-section, but not on its size proper. 
It may take different values for channels with, say, a circular versus a square cross-section, but does not depend on the radius of the circle or the  side length of the square.
\end{itemize}

As will be seen, $h$,  $\lambda_1$ and $\lambda_2$ are obtained relatively easily from an explicit description of the surface micro-relief. $C$
will be given in the form of a power series in Theorem \ref{Ph}. An approximate value for circular channels is $C_{\tiny circle}\approx 0.68$. (See Figure \ref{Ph}.)  Our analysis of $P$ in the weak scattering limit, under the assumption of
isotropic scattering, shows that
$$\vartheta= {h \lambda}/{C}, $$
where we write $\lambda$ for the common value of $\lambda_1$ and $\lambda_2$. See Equation (\ref{O}) and, for a precise mathematical statement, Theorem \ref{Ph} near the end of the paper. We refer to this
equation as the $\vartheta$-{\em factorization formula}.

\subsection{Effective microstructure}
\begin{figure}[h]
\begin{center}
\includegraphics[width=0.7\columnwidth]{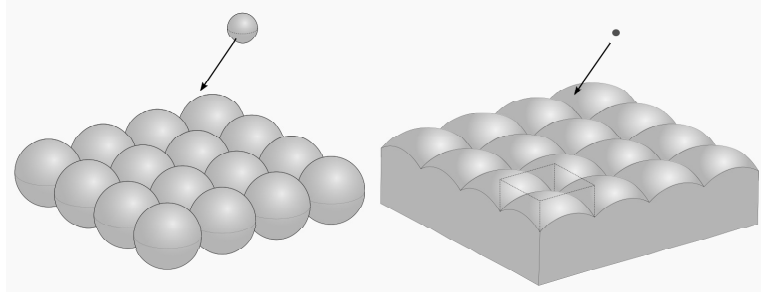}
\caption{  {\small The effective microgeometry of the channel surface depends on the size of the gas molecule. To the   simple surface relief on the left  one associates the effective, or {\em billiard} geometry, on the right. The latter  is the boundary of the region in space that the centers of the spherical gas molecule can occupy.   
In this example, the microgeometry  is periodic;  a {\em tile} or {\em cell}  of this periodic relief is highlighted on the right.}
\label{fig:microgeometry}}
\end{center}
\end{figure}

Before elaborating on this factorization of $\vartheta$,
we  note the natural but important point that   diffusivity is affected by both  the actual surface relief and  the shape and size of the gas molecules. 
Clearly, larger (spherical) molecules are  less affected by the surface irregularities than smaller molecules in the sense that the spread of angles after collisions is wider for the latter. 
 Thus when we speak of the surface microgeometry we have in mind this combined surface-molecule geometry. This  point is illustrated in Figure \ref{fig:microgeometry} for a simple geometry consisting of packed spheres of radius $r_s$ and gas molecules of radius $r_m$.

On the left-hand side of  Figure  \ref{fig:microgeometry}  we have the  picture of a solid surface and a molecule of positive radius impinging on it and on the right the  corresponding effective geometry 
in which the molecule is replaced by a point particle and the surface is thickened accordingly by the amount $r_m$. We call this effective surface geometry the {\em billiard} geometry.
  When probed by spherical molecules of radius $r_m$,  the effective radius of curvature associated to the billiard geometry is $r_s+r_m$. As far as diffusivity is concerned, the same physical  surface should be considered flatter (or less rough) for larger gas molecules. The analysis of the example of Figure \ref{fig:microgeometry} will show, in particular, that the same surface can generate faster or slower gas diffusion depending on the size (larger or smaller) of the gas molecules, at least in the elastic approximation we are considering.   It follows from this discussion that our notion of roughness, as captured by the flatness parameter $h$ to be introduced shortly, is not intrinsic to the physical surface but applies to the combined system. When referring to surface geometry in the subsequent sections, we always have in mind the billiard geometry.

The analytical approach we explore here  applies to much more general situations that   allow  for non-periodic microgeometries and thermal interaction, but for the sake of making the main ideas most transparent we restrict attention to the setting of rigid and periodic structures. This general approach is based on approximating $P-I$, where $I$ is the identity operator,  by a
differential operator $\mathcal{L}$ which we have called in \cite{FNZ2013} {\em MB (Maxwell-Boltzmann)-Laplacian}. 
Under the assumptions of the present paper regarding the surface-molecule interaction, i.e.,  that we have elastic collisions against a rigid surface with no thermal interactions and that diffusivity is {\em isotropic}, the MB-Laplacian turns out to be a generalized Legendre operator in dimension $2$, whose spectral theory is well-understood. As the spectral theory of more general   MB-Laplacians becomes better understood,  it will be possible to undertake a similarly detailed  analysis of  more general situations than discussed in this paper.

\subsection{Definition of the roughness/flatness parameter $h$}
\begin{figure}[h]
\begin{center}
\includegraphics[width=0.4\columnwidth]{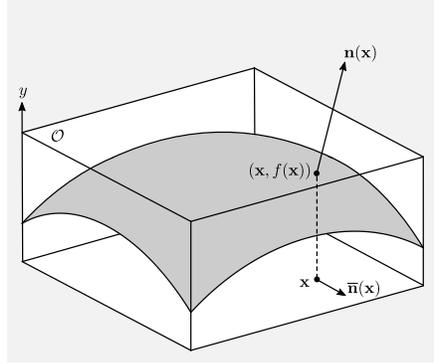}
\caption{{\small For the definition of $h$ and  $\Lambda$ we need the surface function $f$ and the vector field $\overline{\mathbf{n}}$ defined  on $\mathcal{O}$ (or on any rectangle parallel to it). This vector field is  the orthogonal projection to   $\mathcal{O}$ of the unit normal vector field $\mathbf{n}$ on the surface, as indicated in this figure.}}
\label{cell function}
\end{center}
\end{figure}

The precise definition of low roughness will be given  in terms of what we have called in \cite{CFG2021} the {\em flatness parameter} $h$, defined as follows. Let us assume that the surface relief in a billiard cell (for  the example of Figure \ref{fig:microgeometry}, on the right-hand side, this corresponds to a single bump) is represented by the graph of a function $f$ as shown  in Figure \ref{cell function}.  Then $h$ is  
the  maximum value of the square length of the gradient of $f$:  
\begin{equation}\label{h} h=\max_{\mathbf{x}\in \mathcal{O}}|\text{grad}_{\mathbf{x}}f|^2,\end{equation} where the maximum is taken over the points $\mathbf{x}$  in a rectangle parallel to the opening of the cell, denoted by $\mathcal{O}$ in Figure \ref{cell function}.

We should emphasize here the distinction between our flatness coefficient $h$ and standard measures of roughness used today in industry. The  parameter most commonly used to specify surface texture is the {\em average roughness}. It is defined as the average deviation of a surface relief from its mean height. More precisely, let $f_0$ denote the function whose graph gives the physical surface relief, as opposed to $f$, which represents the effective microgeometry. In Figure \ref{fig:microgeometry}, $f_0$ defines the surface on the left-hand side (the part of  the spheres packing surface that is exposed to collisions with gas molecules) and $f$ the one on the right. A commonly used definition of surface average roughness, adapted to our setting of periodic microstructures, is
$$ \text{Ra}=\frac1A \iint_{\mathcal{O}} |f_0(x,y)-\langle f_0\rangle|\, dx\, dy$$
where $\mathcal{O}$ is a rectangular unit of the horizontal plane containing a surface cell, as in Figure \ref{cell function}, $A$ is the area of the rectangle and $\langle f_0\rangle$ is the average height. 
Other measures, such as $\text{Rms}$ that evaluates the root mean square   height profile, are also used.
Our definition of $h$, on the other hand, applies to the effective microgeometry that characterizes the combined surface-molecule system and, as such, it is not directly comparable with standard definitions.
Note that the standard measure of roughness has physical unit of distance while our $h$ is scale-invariant. 

It is possible   to modify the standard definition ($\text{Ra}$) so as to apply to the effective geometry (replacing $f_0$ with  $f$) and perhaps make the resulting parameter scale invariant. 
However, any notion of roughness similar to $\text{Ra}$ would likely not be compatible with our $h$ since  $f$ can have very small height variation between peaks and valleys   and, at the same time, steep gradient over a small area, yielding  large $h$.  
 The parameter $h$ (together with $\lambda$; see the next subsection) may be regarded as a dimensionless measure of the  curvature of the effective geometry. It is surface curvature rather than average height variation  that  is more directly relevant to the determination of diffusivity. 
 
\begin{figure}[h]
\begin{center}
\includegraphics[width=0.6\columnwidth]{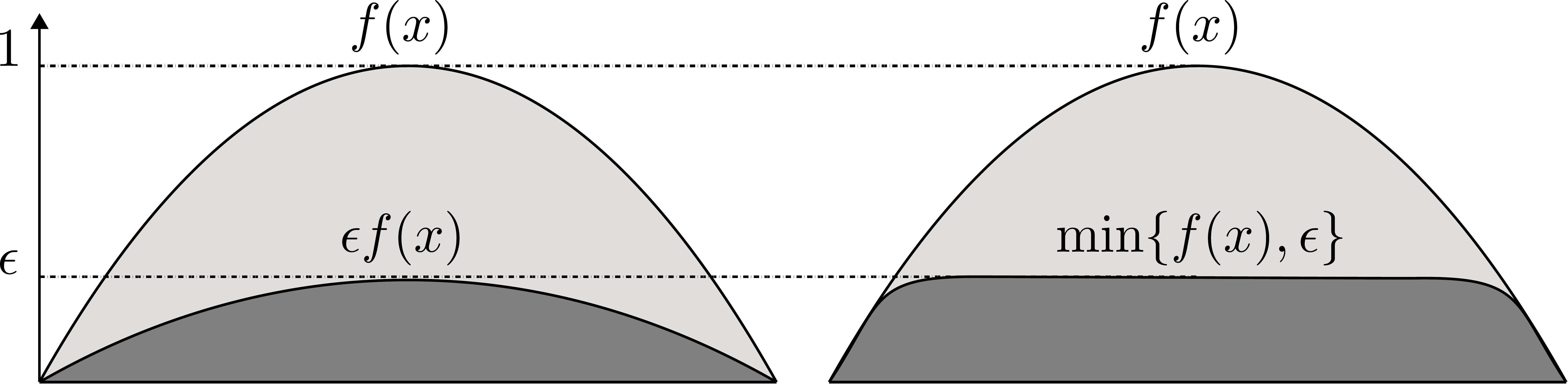}
\caption{{\small Two models of surface polishing. On the left: uniform polishing by a factor $\epsilon$, producing a reduction in the flatness parameter $h$ by $\epsilon^2$. On the right: cut-off at $\epsilon$ does not affect $h$. In both cases, standard measures   such as $\text{Ra}$ give approximately the same reduction in surface roughness.}}
\label{polishing}
\end{center}
\end{figure}

Figure \ref{polishing} may help to clarify the distinction between $h$ and $\text{Ra}$. If the polishing process reduces $\text{Ra}$ without eroding steep slopes, $h$ may remain relatively large without greatly affecting the value of diffusivity. In short, measures of roughness that more directly affect diffusivity must relate to the gradient of $f$ rather than the range of values of $f$. If we control the derivative we also control the value (i.e., small $h$ implies small Ra) but the converse does not necessarily hold. 
The virtue of our flatness parameter $h$ is that it more directly relates geometric properties of the surface and the characteristics of gas transport, at least when $h$ is sufficiently small.  It is an interesting practical issue to investigate what existing  technologies  of surface polishing at the nanoscale can reduce $h$ (and consequently Ra)
rather than just Ra, although it is also a question outside the scope of the present study. 

\subsection{The shape matrix $\Lambda$ and the factorization of $\vartheta$}
For small values of $h$\----corresponding to well-polished, or relatively flat, surfaces\----,
all the geometric parameters characteristic of the surface microstructure having relevance for the diffusion process are contained (in dimension $3$) in $h$ itself and in a real, symmetric $2$ by $2$ matrix $\Lambda$ also obtained from the cell-defining function $f$. (See Figure \ref{cell function}.) In order to define $\Lambda$, let us first introduce a matrix $A$ as follows: for any vector $u$ in the plane tangent to the channel surface at a given point (we should think of $u$ as the tangential component of a post-collision velocity at that point) we set
$$Au := \frac1{\text{Area}(\mathcal{O})}\int_\mathcal{O} \langle \overline{\mathbf{n}}(\mathbf{x}), u\rangle\overline{\mathbf{n}}(\mathbf{x})\, d\mathbf{x}.  $$
Here $\overline{\mathbf{n}}$ is the tangential component of the unit normal vector to the graph of $f$, shown in Figure \ref{cell function},  and
$\langle \overline{\mathbf{n}}(\mathbf{x}), u\rangle$ is ordinary dot product. We then define 
$$\Lambda := \lim_{h\rightarrow 0} A/h $$
when the limit exists. We will call $\Lambda$ the {\em average curvature  matrix} (for small $h$) and its eigenvalues $\lambda_1$, $\lambda_2$ the {\em shape parameters}.

We say that the diffusion is {\em isotropic} if $\Lambda$ is a scalar matrix. Notice that this is a condition on the microstructure and not, naturally, on the transport at the  scale of the channel. Let  $\lambda$ be the single eigenvalue of $\Lambda$.  Then our main remark is that  $\lambda h/C$ provides the appropriate replacement for the Maxwell-Smoluchowski parameter $\vartheta$. 
More precisely, 
our central result  is the $\vartheta$-{\em factorization formula}:
\begin{equation}\label{O}\eta=\frac{2-\vartheta}{\vartheta} +O\left(h^{1/2}\right), \ \ \vartheta = \frac{\lambda h}{C}. \end{equation}
 We describe below how $h$ and $\lambda$ are computed. The value of $C$ is given later in the paper in series expansion form  using generalized Legendre functions.

As an  example, consider for a fixed $f$ the family of microgeometries  $f_\epsilon(\mathbf{x})=\epsilon f(\mathbf{x})$ parametrized by the positive parameter $\epsilon$ and defined over  $\mathcal{O}=[-c_1,c_1]\times [-c_2,c_2]$. 
We may interpret reducing $\epsilon$ as providing a simple mathematical model for the polishing of a surface microgeometry defined by $f$. (See Figure \ref{polishing}.)
For $\epsilon$ small,  
$$ \overline{n}(\mathbf{x})=-\epsilon \text{grad}_{\mathbf{x}} f + O(\epsilon^3), \ \ h=\epsilon^2 \max |\text{grad}_{\mathbf{x}}f|^2.$$
Then the entries of $\Lambda$ are
$$\Lambda_{ij} = \frac{1}{\text{max}|\text{grad}_{\mathbf{x}}f|^2}  \frac1{4 c_1 c_2 } \int_{-c_1}^{c_1}\int_{-c_2}^{c_2} \frac{\partial f}{\partial x_i} 
 \frac{\partial f}{\partial x_j}\, dx_1\, dx_2. $$

For a more explicit example, consider the
model shown in Figure \ref{fig:microgeometry}. 
Let
 $r_s$ be the radius of the spheres  making up the surface and 
$r_m$  the radius of the gas molecules.  
Then, disregarding terms of $4$th order in $r_s/(r_s+r_m)$,  one easily computes
 $$\lambda h\approx \frac13\left(\frac{r_s}{r_s+r_m}\right)^2.$$
 See  Section \ref{example} for details. Thus, insofar as diffusivity is concerned, roughness is lower when the same channel is transporting  gas with larger molecular diameter. 
 
 \subsection{Estimating diffusivity using the factorization equation} \label{physical example} Let us apply   our main result (the $\vartheta$-factorization formula (\ref{O})) to a simple but illustrative example based on the geometry of Figure \ref{fig:microgeometry}, 
 consisting of a planar packing of spheres of radius $r_s$ and spherical gas molecules of radius $r_m$. We assume that the channel has circular cross-section. 
Then $r_s$ and $r_m$ are the only geometric parameters characterizing the microgeometry and $C\approx 0.685$.  The determination of the value of $C$ will be explained in Subsection \ref{C}.
Based on the remarks of the previous subsection,
$$\lambda=\frac16, \ \  h=2\sigma^2+O\left(\sigma^4\right), \  \ \sigma=\frac{r_s}{r_s+r_m}, \ \ C_{\text{\tiny circle}}\approx 0.685$$
so  the Maxwell-Smoluchowski parameter is
$$\vartheta\approx  0.49  \left(\frac{r_s}{r_s+r_m}\right)^2.$$ 
The diffusivity enhancement coefficient $\eta$ then becomes, for low roughness (up to higher orders in $h$),
$$\eta=\frac{2C}{\lambda h}-1\approx 4.10\left(1+\frac{r_m}{r_s}\right)^2-1. $$
This rather simple model can be  a  guide to an important factor determining   $\vartheta$ and $\eta$. Taking $r_s$ as a proxy for the scale of surface irregularities, the squared quotient in this expression is close to $1$ when the size of the gas molecules is significantly smaller than $r_s$; it is approximately $1/4$ when $r_m$ and $r_s$ are comparable in value, and it can be very small when $r_m$ is significantly larger than $r_s$.  Notice that $\eta$ changes by an order of magnitude ($\approx 11.6$) as  $r_m/r_s$ varies from $0$ to $2$.  It is important to note, however, that our approximation is not assured for small values of the ratio $r_m/r_s$ since the flatness parameter $h$ grows larger than  $1$ as this ratio of radii becomes less than approximately $0.4$.  

Taking as a physical example  the diffusion of argon in  carbon nanotubes, we have the following very crude estimate: $r_s$ may be taken to be the radius of carbon, approximately $0.1$ nm,  and $r_m$ the radius of argon molecule, which is approximately $0.18$ nm. In this case $r_m/r_s$ is approximately $1.85$ and $h$ is approximately $0.25$, giving  $\mathcal{D}\approx 32 \mathcal{D}_K.$

Such a large value should be compared with experimental results in \cite{H2006} for airflow through
carbon nanotube membranes, in which flow enhancements $\mathcal{D}/\mathcal{D}_K$ between $16$ and $120$  are observed.   Thus our method of analysis is supported by available data. 
For perspective on what to  expect in more ordinary settings, the tangential momentum accommodation coefficient  (TMAC) which we estimate to be $\vartheta\approx 0.06$ in the present example, is typically close to $1$. In fact, based on a survey of experimental results \cite{AP}, the authors state:
``For monoatomic  gases, we recommend TMAC of $0.926$ for the purpose of applying boundary condition in theoretical analysis and numerical 
computations, for all rare gases, the entire range of Knudsen numbers, and for most surfaces (especially glass, silicon, and steel; exception is platinum).''

We expect that our values for $\vartheta$ and $\eta$ as given by the above formulas are sufficiently accurate to be useful in this   situation, although more detailed work is needed to effectively  implement such a case study.   
 
 \subsection{The method of diffusion approximation in velocity space}
This   topic was developed in greater detail,
  in dimension $2$,  in \cite{CFG2021}. In the present paper, we explain the necessary modifications needed for the more realistic $3$-dimensional context.
The central theoretical ideas   are contained in Theorems \ref{L approx} and \ref{Ph}.  The essential point is that, for relatively weak scattering, the Markov chain generated by $P$, which we now denote by $P_h$ to make explicit the dependence on the flatness parameter,
can be approximated by a diffusion process in velocity space  (not to be confused with the actual gas diffusion in the channel).  The generator of the diffusion process is the MB-Laplacian  $\mathcal{L}$, a differential operator  which is related to $P_h$ according to 
\begin{equation}\label{Leg} P_h\approx I + h\mathcal{L}\end{equation}
in the precise sense of Theorem \ref{L approx}.  One then obtains $\eta$ via the solution of a Markov-Poisson equation by expanding the solution in
terms of the eigenfunctions of $\mathcal{L}$. In the isotropic case, $\mathcal{L}$ has the form
$$\mathcal{L}\Psi(u)= 2\lambda\, \text{div} \left(\left(\rho^2-|u|^2\right) \text{grad}\Psi(u)\right). $$
 This leads to the conclusion of Theorem \ref{Ph}, which is the paper's  main result. The proof of Theorem \ref{Ph} is similar to that of the corresponding result in dimension $2$ from \cite{CFG2021} and is not presented here. The main difference is that in dimension $2$ one is dealing with the standard Legendre (ordinary differential) operator, whereas here the spectral theory of the above partial differential operator is needed instead. 

\begin{figure}[h]
\begin{center}
\includegraphics[width=0.9\columnwidth]{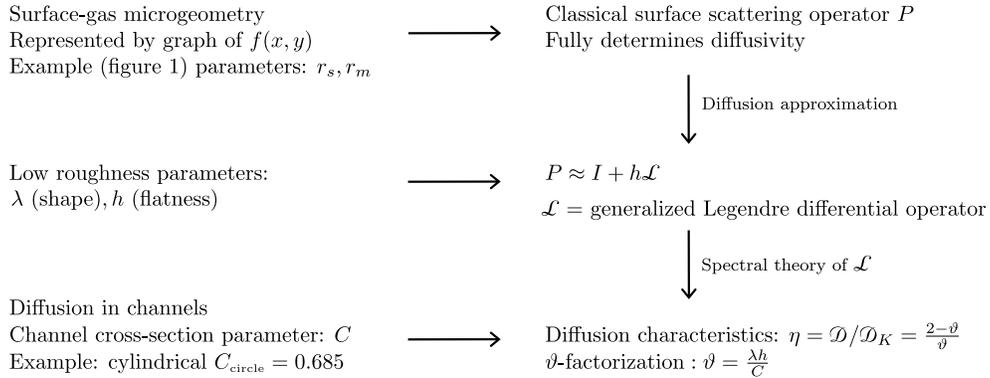}
\caption{{\small Symbols and relations among the main quantities, functions, and operators used in the paper. }}
\label{operators}
\end{center}
\end{figure}

It is important to observe that our use of the generalized Legendre differential operator and its eigenfunctions is not at all arbitrary. It is dictated by the remarkable fact that such operator arises naturally in the diffusion approximation of the Markov operator $P$. Thus, the overall logic in our approach (described diagrammatically in Figure \ref{operators}) can be summarized as follows: as already noted (Equation (\ref{DDK})), $\eta=\mathcal{D}/\mathcal{D}_K$   is fully specified by the spectrum of $P$; for small $h$, $P$ is well approximated by $\mathcal{L}$ up to a multiplicative constant, as in Equation \ref{Leg}. The spectral theory of $\mathcal{L}$ is known explicitly, thus yielding an effective method for obtaining $\eta$. From  this analysis results Equation (\ref{O}). 

When the channel surface microstructure is not rigid but contains moving parts,  $\mathcal{L}$   is a partial differential operator that generates a positive recurrent diffusion in velocity space  whose  stationary probability is the surface Maxwellian (which contains, in particular, Knudsen's cosine law)   at a given temperature   as discussed in detail and in great generality in
\cite{FNZ2013}.  
Our approach to determine  diffusivity for weak surface-molecule interactions  applies to this more general case, except that the spectral theory for the general MB-Laplacian, needed for Theorem \ref{Ph}, is still to be  developed. 

The  diagram in Figure \ref{operators} summarizes  the main notation.

In the following sections, we describe our methods and results in greater mathematical detail.
\section{Random billiards and the  surface operator $P$}
\label{microstructure}
\begin{figure}[h]
\begin{center}
\includegraphics[width=0.4\columnwidth]{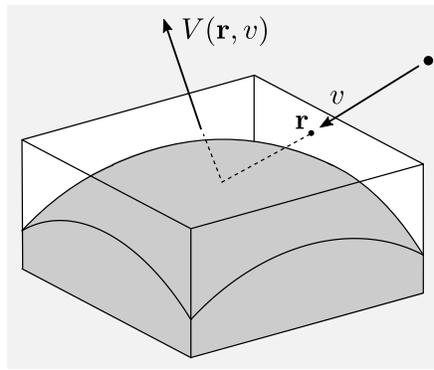}
\caption{  {\small The post-collision velocity $V$ is a function of the pre-collision velocity $v$ and the random variable $\mathbf{r}$ describing the position on the opening of the billiard cell.}
\label{cell}}
\end{center}
\end{figure}

The surface microgeometry and Knudsen diffusivity are  naturally mediated by  a classical scattering operator, here denoted by $P$.
For much of what we show in this paper, we may allow
 the microgeometry to   be very general, although  it will greatly simplify the discussion to assume that it is periodic. The channel surface is then considered to be tiled by identical elements that we call  {\em billiard cells}. A representative cell
is imagined to be contained in a rectangular region, as in Figure \ref{cell}, that has one face open to the interior of the channel (the {\em opening}) and four faces 
on which we impose periodic conditions.

\subsection{Definition of the surface operator}
  Since we regard these cells to be very small relative to  the diameter of the channel,  it makes sense to characterize a molecule-surface collision as a scattering event given as follows: when a gas molecule impinges on the surface,
its post-collision velocity $V(\mathbf{r},v)$  is a function of the incoming velocity $v$
and of a point $\mathbf{r}$ on the opening of the cell which we regard as being random, uniformly distributed on the cell opening. (See \cite{FY2004} for more details.) Let us denote the opening by $\mathcal{O}$. The event that the post-collision velocity lies in a set $S$ of $3$-dimensional vectors will then have the probability
$$\text{Prob}(V \text{ lies in } S\  | \text{ pre-collision velocity is }  v) = \frac{1}{\text{Area}(\mathcal{O})} \iint_\mathcal{O} \mathbbm{1}_S(V(\mathbf{r},v))\, d\mathbf{r} =\mathbb{E}_v[\mathbbm{1}_S(V)].$$ 
Here $\mathbbm{1}_S$ is the {\em indicator function} of $S$ which, by definition,  takes on the value $1$ in $S$ and $0$ outside $S$, and
$\mathbb{E}_v$ indicates conditional expectation given the pre-collision velocity $v$. Using more general test functions $f(V)$, we define the classical scattering operator $P$ (or Markov transition operator) associated with the  microgeometry as
\begin{equation}\label{P}(Pf)(v)= \mathbb{E}_v[f(V)]. \end{equation}

All the  geometric properties of the surface that are relevant to diffusion are contained in $P$. It is the surface's scattering signature.   Given $\mathbf{r}$ and $v$, $V(\mathbf{r},v)$ is the result of the deterministic motion involving one or more billiard-like, specular collisions in the billiard cell. The random nature of $V$ is due to the randomness of the point of incidence $\mathbf{r}$ in the opening of the billiard cell. Notice how the random flight in the cylindrical channel is then  obtained in this surface model by combining two essentially independent steps: 
one step is the  free flight between two consecutive molecule-surface collisions and the other is the determination of the post-collision velocity as a random function of the pre-collision velocity. The latter step does not  require knowing the exact position at which the molecule hits the surface.

\subsection{General properties of the surface operator}

As has been described elsewhere (for example, \cite{CFZ2016,CF2012,Feres2007}), the operator $P$ has many nice properties, a few of which we summarize here:

\begin{itemize}
\item {\em Knudsen cosine law.}  The probability distribution $$d\mu=\frac1{2\pi}\cos \varphi \, d\Omega$$ on solid angles, where $\varphi$ is the angle a post-collision velocity makes with the normal vector to the channel surface, is stationary under the Markov chain defined by $P$, irrespective of the given microgeometry.
\item {\em Self-adjointness of $P$}. On the Hilbert space of  functions  $f(V)$ that are square integrable with respect to the Knudsen cosine distribution $\mu$, $P$ is a bounded self-adjoint operator. Its spectrum is real and is contained in the closed interval $[-1,1]$.
\item {\em Spectrum.} The spectrum of $P$ (real and often discrete) is thus a signature of the microgeometry.  It   fully determines diffusivity. (See \cite{FZ2012} for many examples. See also \cite{CFZ2016}, where the  integral formula (\ref{DDK}) over the spectrum of $P$ is obtained for the diffusivity.) Of special importance is the {\em spectral gap} $\gamma$, defined as the gap between the top eigenvalue $1$ and the rest of the spectrum of $P$. This quantity has a pronounced effect on diffusivity, as will be further noted in this paper.
\end{itemize}

This model of microstructure can be extended to allow for moving parts and surface potentials. In this more general situation, the stationary probability distribution
on velocity space is the (boundary) Maxwell-Boltzmann distribution at a given temperature, expressed in terms of the variance of
velocities of the moving parts of the surface microstructure. This  contains the Knudsen cosine law as the distribution of scattering directions. (See \cite{CFZ2016,CF2012,FNZ2013}.)

 \section{Random flight in a cylindrical channel}
 We summarize in this section some results about diffusion,  mostly  special cases of theorems proved in \cite{CFZ2016}, that are needed for the present paper. We assume for concreteness that the channel cross-section is circular. Different cross-sections correspond to different displacement functions, denoted $X$ below. All the other elements of the analysis remain unchanged.
 
  \begin{figure}[h]
\begin{center}
\includegraphics[width=0.9\columnwidth]{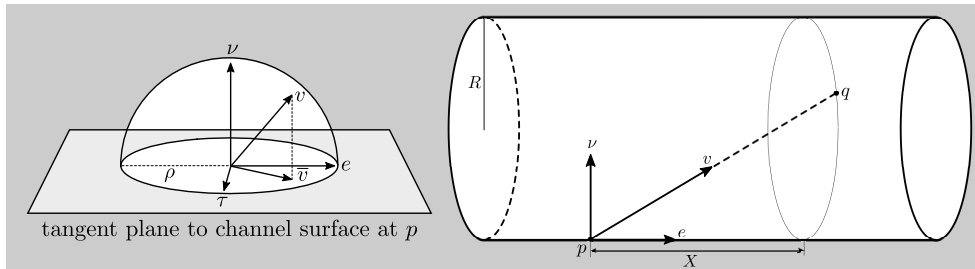}
\caption{{\small On the left is the definition of the moving frame $(\tau, e, \nu)$ on the channel surface.  Pre- and post-collision velocities, $v$ and $V$, lie on (the surface of) a sphere of radius $\rho=|v|$. It is convenient to represent   velocities  by their orthogonal projections $\overline{v}$ and $\overline{V}$ to the disc $D_\rho$ of radius $\rho$ perpendicular to  $\nu$. On the right is the definition of the displacement function $X=(q-p)\cdot e$ where $q=p+tv$ is the next point of collision, $t$ is the time between collisions, and $u$ is the projection of $v$ to the disc of radius $\rho$ perpendicular to $\nu$.}}
\label{hemisphere}
\end{center}
\end{figure}

 \subsection{The Markov chain of scattered velocities and random flight}
 The  remarks of Section \ref{P} can be expressed in the language of Markov chains. 
Let $V_1, V_2, V_3, \dots$ be the sequence of postcollision velocities of a tagged gas molecule as it undergoes  random flight in the cylindrical channel of radius $R$. This sequence of random variables constitutes a Markov chain with transition operator $P$. In the stationary regime, even though  $V_i$  and $V_{i+1}$ can be strongly correlated, these random variables satisfy the cosine distribution of directions regardless of the given microstructure. Given $V_i$, the distribution of $V_{i+j}$ converges to the cosine law at an exponential rate as $j$ increases. This rate of relaxation, which is  dominated by the spectral gap $\gamma$, is a key factor influencing diffusivity.    
 The mathematical details in dimension $2$ are developed in \cite{CFG2021}.  
 
 As the speed $\rho=|V_i|$ does not change during the process due to the assumption that collisions are elastic, $V_i$ is determined by its orthogonal projection  $\overline{V}_i$ to the disc $D_\rho$ of radius $\rho$: $$V_i=\left(\overline{V}_i, \sqrt{\rho^2-|\overline{V}_i|^2}\right).$$ It is not difficult to show that having the cosine law for the stationary probability distribution of the process $V_i$ is equivalent   to  the process $\overline{V}_i$ having  the uniform distribution on $D_\rho$ for its stationary probability  distribution.

 We introduce the orthonormal moving frame $(\tau, e, \nu)$ on the surface of the cylindrical channel defined  in Figure \ref{hemisphere}. At any given point, $\tau$ is the unit vector tangent to the circle cross-section, $e$ is the unit length vector pointing along the axis of the cylinder, and
 $\nu$ is the unit vector normal to the cylinder surface,  pointing in. 
If $P_1, P_2, \dots $ are  the collision points  with the channel surface of a random flight trajectory, where $i$ indicates the flight step, we write $$\tau_i=\tau(P_i),\ \  e_i=e(P_i)=e \ (\text{constant}),\ \  \nu_i=\nu(P_i).$$

We now describe the mathematical model of the random flight in a cylindrical channel.  Each step of the random flight is determined
 as follows.
 \begin{itemize}
 \item  Let $P_i$ be the point on the channel wall at the moment $t_i$, where $i$ indicates the present collision step of the random flight.
 At this step, let $V_i$ be the post-collision velocity. 
 \item The time $t_{i+1}= t_i + T_i$ of the next collision  is obtained from 
 $$ T_i = 2R\frac{V_i\cdot \nu_i}{(V_i\cdot \nu_i)^2+(V_i\cdot \tau_i)^2}.$$
 The next collision point is $P_{i+1}=P_i+T_i V_i$ and the displacement along the axis of the channel accounted for by this segment of random flight is
 $X_i =  (V_i\cdot e) T_i$.
 \item Finally, the post-collision velocity $V_{i+1}$ is  obtained  as follows: Given $V_i$ and a random point $\mathbf{r}$ in $\mathcal{O}$ uniformly distributed, we obtain $V_{i+1}=V(\mathbf{r},V_i)$, as indicated in Section \ref{microstructure}, for a given choice of microstructure. 
 \item The total displacement after $n$ steps of the random flight is then
 $S_n =X_0+\cdots + X_{n-1}. $
 \end{itemize}

Introducing the  displacement  $X(p,v)$  and  the free-flight time $T(p,v)$ functions   
   $$ X(p,v)=2R\frac{(v\cdot \nu) (v\cdot e)}{(v\cdot \nu)^2+(v\cdot\tau)^2}, \ \ T(p,v)=2R\frac{v\cdot \nu}{(v\cdot \nu)^2+(v\cdot\tau)^2},$$
where $\tau, e, \nu$ are evaluated at $p$, lets us write $X_i=X(P_i, V_i)$ and $T_i=T(P_i,V_i)$. 
It is clear that, relative to the stationary probability distribution, $X$ has mean $0$. A simple integral calculation gives the variance of $X$ and mean value of $T$:
$$\mathbb{E}\left[X^2\right]= \frac{4}{3} R^2, \ \ \mathbb{E}[T]=\frac{2R}{\rho}. $$

\subsection{Channel diffusion as   limit of the random flight process}
The diffusion process arises from the random flight model through an application of a Central Limit Theorem. In this regard, the facts we need for this paper will be taken from \cite{CFZ2016}. A concrete way to think about the CLT in the present context comes by the consideration of the following idealized experiment. Keeping in mind that one quantity that can be effectively measured
is the time that molecules undergoing random flight take to escape from a finite length channel,  let the channel be  a cylinder of radius $R$ and total length $2L$, and let 
$\tau(R,L,\rho)$ be the mean exit time (where the random flight is defined for some choice of microstructure with Markov operator $P$), assuming starting point at the middle of the channel and random initial velocity (with speed $\rho$). The CLT (as in 
\cite{CFZ2016}) implies the asymptotic expression for large $L$:
$$\tau(R,L,\rho) \sim \frac{L^2}{\mathcal{D}}  $$
where $\mathcal{D}$ is the constant of self-diffusivity.  A simple dimensional analysis argument shows that $\mathcal{D}=C_P R\rho$, where
$C_P$ is a constant, depending only on the microstructure, obtained by taking the limit for large $a$ of  the dimensionless quantity
$a^2/F(a)$, where $F(L/R)=(\rho/R)\tau(L,R,\rho)$.

A standard way to express $\mathcal{D}$  is through the relation $\mathcal{D}=\eta \mathcal{D}_K$, where $\mathcal{D}_K$  ($K$ standing for {\em Knudsen}) is the diffusivity obtained under the assumption that the velocity process is independent and  identically (cosine law) distributed. One finds by an application of the standard CLT:
\begin{equation}\label{standard}\mathcal{D}_K=\frac{2}{3}R\rho. \end{equation}

We give now a formula for $\mathcal{D}$ for a given choice of $P$.  As already noted, $P$ is a self-adjoint operator on the Hilbert space $\mathcal{H}$ of square integrable (complex-valued) functions on $D_\rho$ with inner product
$$\langle f, g\rangle =\int_{D_\rho} \overline{f(u)} g(u) du. $$
(Recall that the cosine law on the hemisphere of radius $\rho$ corresponds to the uniform probability distribution on $D_\rho$.) 
The norm derived from this inner product will be denoted $\|f\|=\sqrt{\langle f,f\rangle}$.
The spectral theorem for self-adjoint operators provides a projection-valued measure $\Pi(d\lambda)$ on the spectrum of $P$, a closed subset of the real line. Combining $\Pi$ with the displacement function $X$ (which is square-integrable due to its finite variance; here, the fact that the channel is a cylinder is being used; this would not be the case for diffusion between two parallel plates)  one obtains a finite measure on the spectrum:
$$\Pi_X(d\lambda) =\|X\|^{-2}\langle X, \Pi(d\lambda)X\rangle. $$
We can now state a fundamental formula showing that the information about the microstructure of the channel surface relevant to 
diffusion properties is contained in the spectrum of $P$:
\begin{equation}\label{spectral formula} \eta=\int_{-1}^1 \frac{1+\lambda}{1-\lambda}\, \Pi_X(d\lambda),\end{equation}
recalling that $\mathcal{D}=\eta\mathcal{D}_K$. (This is a reformulation of the integral in  Equation (\ref{DDK}).) The argument used in \cite{CFZ2016} to derive (\ref{spectral formula})   actually shows
the following:
\begin{equation}\label{key formula}
\eta =1+2 \|X\|^{-2}  \left\langle X, P(I-P)^{-1}X\right\rangle.
\end{equation}
Note that $\|X\|=(2/\sqrt{3})R$ as seen above. Equation (\ref{key formula}) together with a canonical approximation of $P$ by a generalized Legendre differential operator, discussed in the next section, are the key ingredients for the main results of the present paper.

To understand how this equation can be used to find $\eta$, let us introduce the function $Y=(I-P)^{-1}X$, which is a solution 
of the Markov-Poisson equation 
$$ LY = -X,$$
in which $L=P-I$ may be called the  {\em Laplacian} associated to $P$.  As will be seen next, when the channel surface has low roughness (i.e., it is a well-polished surface, defined mathematically by the small value of the flatness parameter $h$), $L$ is  well-approximated by a differential operator $\mathcal{L}$ for which the solution to the  equation $\mathcal{L}Y=-X$ can be effectively obtained.  Here $\mathcal{L}$ is a generalized Legendre operator, which we describe in the next section. 
  
 \section{The shape matrix $\Lambda$ and the generalized Legendre operator $\mathcal{L}$}\label{example}
 We specialize here a key remark made in much greater generality in  \cite{FNZ2013}.
 When surface roughness is small (well-polished surfaces), it will be seen that all the micro-geometric parameters influencing diffusivity, assuming isotropic diffusion,
 are summarized  by $h$ and a matrix $\Lambda$, which we call the {\em shape matrix}. For sufficiently flat billiard microgeometries we may assume that the surface in a billiard cell is   the graph of a function $f$, that is,    points
 $\mathbf{r}$ on the surface take the form $\mathbf{r}=(\mathbf{x},f(\mathbf{x}))$, where $\mathbf{x}$ may be regarded as a point in  the cell opening $\mathcal{O}$ (or, in fact, any vertical translate of $\mathcal{O}$; see Figure \ref{cell function}). 
 For any given $\mathbf{x}$ let $\mathbf{n}(\mathbf{x})$ be the unit length perpendicular vector to the surface at $(\mathbf{x},f(\mathbf{x}))$ and let $\overline{\mathbf{n}}(\mathbf{x})$ be the orthogonal projection of $\mathbf{n}(\mathbf{x})$ to the plane of $\mathcal{O}$. 
 Since $\mathbf{n}(\mathbf{x})$ lies in the direction of the gradient of $F(\mathbf{x},y)=y-f(\mathbf{x})$, and denoting by $e=(0,0,1)$ the unit vector in the direction of the $y$-axis (see Figure \ref{cell function}), we have
 $$\mathbf{n}(\mathbf{x}) =  \frac{e-\text{grad}_\mathbf{x} f}{\sqrt{1+|\text{grad}_{\mathbf{x}}f|^2}}, \ \ \ \overline{\mathbf{n}}(\mathbf{x}) =  -\frac{\text{grad}_\mathbf{x} f}{\sqrt{1+|\text{grad}_{\mathbf{x}}f|^2}}.$$

Recall that the   maximum value of $|\text{grad}_{\mathbf{x}}f|^2$ over the points $\mathbf{x}$ in $\mathcal{O}$ is the flatness parameter $h$ of the surface defined by $f$.  
Also recall that  $\Lambda$ is obtained from  the  $2$ by $2$ matrix $A$ via the limit $\Lambda := \lim_{h\rightarrow 0} A/h$, and that $A$
acts on vectors $u$ in the plane perpendicular to $e=(0,0,1)$ as
$$Au := \frac1{\text{Area}(\mathcal{O})}\int_\mathcal{O} \langle \overline{\mathbf{n}}(\mathbf{x}), u\rangle\overline{\mathbf{n}}(\mathbf{x})\, d\mathbf{x}.  $$
The flatter the surface the smaller is $Au$, in such a way   that, in typical cases, we should expect $Au$ to be of the order of the flatness parameter $h$, justifying the definition of $\Lambda$.

 \subsection{Example of computation of microparameters $\lambda$ and $h$}
Let us consider a couple of representative examples. We write $\mathbf{x}=(x_1,x_2)$. Suppose $\mathcal{O}$ is the rectangle    such that $|x_1|\leq c_1$ and
$|x_2|\leq c_2$ for   positive numbers 
$c_1, c_2$. Let $a_1, a_2, b$ be also positive, $\epsilon$ a small number, and set
$$f_\epsilon(\mathbf{x})=\frac{b}{\epsilon}\left\{\sqrt{1-\left(\frac{\epsilon x_1}{a_1}\right)^2  -\left(\frac{\epsilon x_2}{a_2}\right)^2} -
\sqrt{1-\left(\frac{\epsilon c_1}{a_1}\right)^2  -\left(\frac{\epsilon c_2}{a_2}\right)^2} \right\}. $$
The graph of $f_\epsilon$ is a piece of ellipsoid over $\mathcal{O}$ whose principal axes are proportional to $1/\epsilon$.   
Straightforward calculation gives
$$h = b^2\left(\frac{c_1^2}{a_1^4} +\frac{c_2^2}{a_2^4}\right)\epsilon^2 + O(\epsilon^4)$$
and 
$$A/h =\frac13 \left(\frac{c_1^2}{a_1^4} +\frac{c_2^2}{a_2^4}\right)^{-1}\left(\begin{array}{cc}\frac{c_1^2}{a_1^4} & 0 \\0 & \frac{c_2^2}{a_2^4}\end{array}\right) + O\left(\epsilon^2\right).$$
Therefore
$$\Lambda =\left(\begin{array}{cc}\lambda_1 & 0 \\0 & \lambda_2\end{array}\right) \text{ where } \lambda_i = \frac13\frac{c_i^2}{a_i^4} \left(\frac{c_1^2}{a_1^4} +\frac{c_2^2}{a_2^4}\right)^{-1} \text{ for } i=1,2. $$
We record  for later use the isotropic case, in which $$a:=a_1=a_2, \ \ c:=c_1=c_2, \ \ \lambda_1=\lambda_2=1/6, \ \  h=2\left(\frac{b}{a}\right)^2\left(\frac{c}{a}\right)^2 \epsilon^2 + O(\epsilon^4).$$  

For the special case of the example of Figure \ref{fig:microgeometry}, let $r_s$ be the radius of the spheres that constitute the surface and 
$r_m$ the radius of the gas molecules. Then  $\epsilon = r_s/(r_s+r_m)$ and  $a_1=a_2=b=c_1=c_2=r_s$ and we have
$$\lambda h\approx \frac13\left(\frac{r_s}{r_s+r_m}\right)^2,$$
disregarding terms of $4$th order in $r_s/(r_s+r_m)$.

\subsection{The   generalized Legendre operator}
The interest in $\Lambda$ is that it appears as a matrix of coefficients for a differential operator $\mathcal{L}$ that we now define, which generalizes the well-known Legendre differential  operator in dimension $1$. This generalized Legendre operator will appear in the approximation of the Markov operator $P$ for  small $h$ microstructures.
 $\mathcal{L}$ will act on functions  of the molecular velocity $v$.  Since molecular speed, $\rho$, does not change during the random flight inside the channel under the assumption that collisions are elastic, we may identify the space of velocities at any point on the surface of the channel with the sphere of radius $\rho$.  Rather than use $v$ itself, it is convenient to represent velocities by their orthogonal projection to the disc $D_\rho$ of radius $\rho$ perpendicular to the unit normal vector $\nu$ to the surface of the channel as indicated in Figure \ref{hemisphere}.  Representing the orthogonal projection of a velocity $v$ by $\overline{v}$, we wish to define a differential operator on functions $\Psi(\overline{v})$ for $\overline{v}$ in $D_\rho$.   

We can now define
$$(\mathcal{L}\Psi)(\overline{v})= -4 \langle \text{grad}_{\overline{v}}\Psi, \Lambda \overline{v}\rangle + 2\left(\rho -\left|\overline{v}\right|^2\right)\text{Tr}\left(\Lambda \text{Hess}_{\overline{v}}\Psi\right).$$
Here $\text{grad}$ and $\text{Hess}$ are the two-dimensional  gradient and Hessian.  A simpler expression results by using coordinates on the plane adapted to an orthonormal basis of eigenvectors of the symmetric matrix $\Lambda$. Writing   $u=\overline{v}=u_1e_1+u_2 e_2$, where $\Lambda e_i = \lambda_i e_i$,  then one easily shows that
$$\frac12(\mathcal{L}\Psi)(u)=  \lambda_1 \frac{\partial}{\partial u_1}\left(\left(\rho^2-|u|^2\right)\frac{\partial \Psi}{\partial u_1}\right)+
 \lambda_2 \frac{\partial}{\partial u_2}\left(\left(\rho^2-|u|^2\right)\frac{\partial \Psi}{\partial u_2}\right).$$
 
 In this paper we restrict attention to isotropic diffusion, mainly because an explicit spectral theory of the general $\mathcal{L}$ does not seem to be available to the best of our knowledge. This  amounts to assuming that $\Lambda$ is a scalar matrix of the form $\lambda I$ where $I$ is the identity. In this case  
 $$\left(\mathcal{A}\Psi\right)(u)  :=\frac1{2\lambda}\left(\mathcal{L}\Psi\right)(u)= \text{div} \left(\left(\rho^2-|u|^2\right) \text{grad}\Psi(u)\right).$$
Note that the only parameters associated to the billiard microstructure are then $\lambda$ and $h$. Without loss of generality, we set $\rho=1$. The disc of radius $1$ will be written $D$ instead of $D_1$.
The  spectral theory of  $\mathcal{A}$,  just as that of its   one-dimensional counterpart,  is available although not widely  known. See \cite{M1990,S1969}. We summarize  the main facts in the next proposition.

\begin{prop}[Spetral theory of $\mathcal{A}$]\label{propoA}
The eigenvalues of $\mathcal{A}$ are given by
$$ \lambda_{k\ell}=(2\ell+1)(2\ell+2k+1), \ \ k, \ell =0,1, 2, \dots$$
where the multiplicity of $\lambda_{k\ell}$ is $1$ when $k=0$ and $2$ when $k\geq 1$. For each $k,\ell$ and $j=\pm 1$, the corresponding eigenfunctions are given by
$$ \phi_{k\ell j}(u)=F(-\ell, \ell+k+1;k+1;|u|^2)H_{kj}(u),$$
where $F$ is the classical hypergeometric function and $H_{kj}(u)= |u|^k e^{ijk \theta}$, with $\theta$ being the polar angle  of $u\in D$.
 \end{prop}
 Note that $\{H_{kj}\}$ is the basis of the space of harmonic homogenous polynomials of degree $k$. The explicit form of $F$ (traditionally written 
$\tensor*[_2]{F}{}$) is
$$F(-\ell, \ell+k+1;k+1;s) =\sum_{n=0}^l (-1)^n \frac{\binom{\ell+k+n}{n} \binom{\ell}{n}}{\binom{k+n}{n}}s^n.$$

 \section{Differential approximation of the shape operator $P$ and self-diffusivity}
 We can now state the key fact we need that relates the Markov operator $P$ and the differential operator $\mathcal{L}$.  The result quoted here was established in much greater generality in \cite{FNZ2013}, although the observation that this relation can be used to approximate  Knudsen self-diffusivity was not made there. 
 
\subsection{The operator approximation result}
 \begin{thm}\label{L approx}
 Let $(f_h)$ be a family of piecewise smooth functions defined on $\mathcal{O}$ with associated average curvature matrix $\Lambda$,  where $h>0$ is the flatness parameter.  Let $(P_h)$ be the corresponding family of Markov transition operators. Then for any function $\Psi$ on $D$ having continuous derivatives to order at least $3$,
 $$P_h\Psi(u)-\Psi(u) = h \mathcal{L}\Psi(u) + O\left(h^{3/2}\right)$$
 holds for each $u$ such that every initial condition  with velocity $v$ having projection $\overline{v}=u$ results in a trajectory  that collides only once
 with the boundary surface of the cell. Here $\mathcal{L}$ is the generalized Legendre operator associated to $\Lambda$.
 \end{thm}

Using the basis  $\{\phi_{klj}\}$ of eigenfunctions for $\mathcal{L}$, we are able to construct solutions
of the Markov-Poisson equation  $(P_h-I)Y=X$. In analogy with Theorem 4 of \cite{CFG2021} we are thus able to give an approximation formula for $\eta$ in the present $3$-dimensional setting based on  Equation (\ref{key formula}). Roughly, we have
$\eta\approx 1 + 2\|X\|^{-2} \langle X, PY\rangle$ where $Y$ is now solution to $h\mathcal{L}Y=-X$ or
$$2\lambda h \mathcal{A}Y=-X. $$

This analysis leads to the following $3$-dimensional counterpart (proved similarly) to Theorem 8 of \cite{CFG2021}. It is the main result of the present paper.

\subsection{The factorization of the Maxwell-Smoluchowski parameter $\vartheta$}\label{C}
\begin{thm}\label{Ph}
Let $(P_h)_{h>0}$ be a family of random  billiard transition operators for a family of billiard cells satisfying  the geometric assumptions of Theorem \ref{L approx}.  Then 
$$\eta = \frac{2-\vartheta}{\vartheta} + O\left(h^{1/2}\right)$$
where $\vartheta = \lambda h/C$ and  
\begin{equation}\label{tripleS}C=\frac12\sum_{j=\pm1}\sum_{k,\ell=0}^\infty \frac{\left\langle \phi_{k\ell j},\overline{X}\right\rangle^2}{\lambda_{k\ell}\|\phi_{k\ell j}\|^2}\end{equation}
only depends on the macrogeometry of the cylinder channel. (In fact, only on the shape of the cross-section of the channel, up to scale. Thus, for a circular cross-section of radius $R$, $C$ is simply a number independent of $R$.) Here $\overline{X}=X/\|X\|$.
\end{thm}

 \begin{figure}[h]
\begin{center}
\includegraphics[width=0.5\columnwidth]{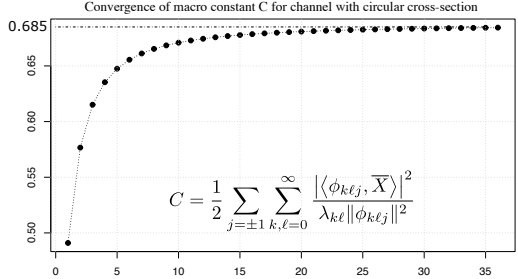}
\caption{{\small Approximation of the macroscopic constant $C$ for a circular straight channel. The   value $C\approx 0.685$ was obtained by truncating the triple iterated infinite series 
for $k,\ell\leq 35$.  More precisely, for each $j\in \{-1,+1\}$, $k, \ell\in \{0, 1, \dots, 35\}$ the terms $A_{jk\ell}$ in the sum are calculated and for each $\ell$ the sum $B_\ell$ of the $A_{jk\ell}$ over $j$ and $k$ are obtained. The plot shows the values of $\sum_{\ell=1}^nB_\ell$ for $n=1,\ldots, 35$. Further details about the calculation of $A_{jk\ell}$ are given in the body of this article. See Theorem \ref{Ph}.}}
\label{convergence}
\end{center}
\end{figure}

    A few words are in order regarding the determination of $C$. The terms of the infinite sum in Equation (\ref{tripleS}) involve the eigenfunctions $\phi_{k\ell j}$ and eigenvalues $\lambda_{k\ell}$ of the operator $\mathcal{A}$ (related to $\mathcal{L}$) given in Proposition \ref{propoA}. The inner product $\langle\cdot, \cdot\rangle$ and norm $\|\cdot\|$ in Equation (\ref{tripleS}) are the standard operations in the Hilbert space  
of square integrable functions on the disc $D_\rho$ indicated in Figure \ref{hemisphere}.   $X$ is the displacement function between two consecutive scattering events of the random flight inside the cylindrical channel, as shown in the same figure. Thus the quantities $\|X\|$, $\|\phi_{k\ell j}\|$ and $\left\langle \phi_{k\ell j},\overline{X}\right\rangle$ require evaluating integrals over $D_\rho$. We used MATLAB for the approximate evaluation of $C_\text{\tiny circle}$ as indicated in Figure \ref{convergence}.

\section{Conclusions} 
We  revisit and refine the classical Maxwell-Smoluchowski theory of gas-diffusion in channels.  
 This refinement consists of a new method for the computation of the tangential momentum accommodation coefficient $\vartheta$ of self-diffusivity of gases in straight channels, in the  large Knudsen number  regime, based on an explicit description 
of the channel surface microgeometry. It is assumed that the molecule-surface interaction is modeled by elastic   collisions of hard spheres against  a rigid surface.
The main conclusion is that, for surfaces having a high degree of polish, or low roughness, as measured by a geometric parameter $h$ which we have called {\em flatness}, the diffusivity enhancement factor  $\eta=\mathcal{D}/\mathcal{D}_K$ can be expressed  as follows: $\eta=(2-\vartheta)/\vartheta$ and $\vartheta=\lambda h/C$. All the  quantities involved   can be calculated from first principles given an explicit  mathematical description of the  surface microgeometry and the shape of the cross-section of the straight channel. Here $C$ is a scale-independent constant that depends only on this  cross-section. For example, for  a circular channel this number  is approximately $0.68$ and is independent of the channel diameter. The flatness parameter $h$, as already noted, gives the overall level of surface polish so that low values of $h$ imply a low degree of surface `roughness'\---- a widely used term that is given precise mathematical meaning in our work through the analytic definition of $h$. And the shape parameter $\lambda$ is a measure of  surface curvature independent of $h$. Both $\lambda$ and $h$ are easily obtained from the mathematical  model of the surface microgeometry whereas $C$ is given by an infinite series as indicated in Theorem \ref{Ph}.  

Applying our main result (Equation (\ref{O}))  to a very simple model of microgeometry, and
taking as a physical example  gas diffusion of argon in  carbon nanotubes, we obtain the following very crude estimate: $r_s$ may be taken to be the radius of carbon, approximately $0.1$ nm,  and $r_m$ the radius of argon molecule, which is approximately $0.18$ nm, giving  $\mathcal{D}\approx 32 \mathcal{D}_K.$  (See the details in Subsection \ref{physical example}.)This large value should be compared with experimental results in \cite{H2006} for airflow through
carbon nanotube membranes, in which flow enhancements $\eta$ between $16$ and $120$  are observed.   
 
 Our method for obtaining $\vartheta$  is based on an approximation of the classical scattering operator, $P$,  that represents the gas-surface interaction by a diffusion operator in velocity space, which we call Maxwell-Boltzmann (MB) Laplacian. In the present work, the MB-Laplacian is a generalized Legendre differential operator. This approximation method  has much greater validity than demonstrated here. For more refined models of gas-surface interactions that allow for energy exchange, the associated MB-Laplacian is also known due to  our earlier work \cite{FNZ2013}, but their spectral theory is not presently well understood. Further progress on this mathematical topic will allow for  greatly extending the applicability of the analysis developed here.

\end{document}